\documentclass[11pt]{article}

\usepackage{amsmath,amsfonts,amssymb,amsthm}
\usepackage[english]{babel}

\newcommand{\sA}{{\mathcal A}}
\newcommand{\sB}{{\mathcal B}}
\newcommand{\sD}{{\mathcal D}}
\newcommand{\sG}{{\mathcal G}}
\newcommand{\sL}{{\mathcal L}}
\newcommand{\sR}{{\mathcal R}}

\newtheorem{theorem}{Theorem}[section]
\newtheorem{lemma}[theorem]{Lemma}

\newtheorem{remark}[theorem]{Remark}

\newtheorem{assumptions}[theorem]{Assumptions}

\newtheorem{definition}[theorem]{Definition}

\numberwithin{equation}{section}

\begin{document}


\centerline{{\Large Optimal boundary control with critical penalization}
\vspace{2mm}} 
\centerline{{\Large for a PDE model of fluid-solid interactions}\footnote{
This research was started while the second author was visiting the
{\em Centro di Ricerca Matematica Ennio De Giorgi} of the Scuola Normale
Superiore in Pisa, whose support is acknowledged.
\\
The first author acknowledges partial support of the Universit\`a degli 
Studi di Firenze, within the 2008 Project 
``Calcolo delle variazioni e teoria del controllo'', 
as well as partial support of the Italian MIUR, under the (PRIN~2007) Project 
``Regolarit\`a, propriet\`a qualitative e controllo di soluzioni 
di equazioni alle derivate parziali non lineari''.
\\
The second author acknowledges support of the National Science Foundation under Grant DMS-060666882.}
\vspace{5mm}} 


\centerline{
{\large Francesca Bucci$^{\textrm{a}}$, Irena Lasiecka$^{\textrm{b}}$}}
\vspace{8mm}      

\noindent
{\small $^{\textrm{a}}$Universit\`a degli Studi di Firenze, Dipartimento di
Matematica Applicata, Via S.~Marta 3, 50139 Firenze, Italy; e-mail
{\tt francesca.bucci@unifi.it}} 
\vspace{2mm}\\  
{\small $^{\textrm{b}}$University of Virginia, Department of
Mathematics, P.O.~Box 400137, Charlottesville, VA 22904-4137, U.S.A.,
e-mail {\tt il2v@virginia.edu}} \vspace{8mm}

\noindent
{\small {\bf Abstract:}
We study the finite-horizon optimal control problem with quadratic functionals
for an established fluid-structure interaction model.
The coupled PDE system under investigation comprises a parabolic (the fluid) and a 
hyperbolic (the solid) dynamics; the coupling occurs at the interface between 
the regions occupied by the fluid and the solid.
We establish several trace regularity results for the fluid component
of the system, which are then applied to show well-posedness of the Differential 
Riccati Equations arising in the optimization problem.
This yields the feedback synthesis of the unique optimal control, 
under a very weak constraint on the observation operator; in particular, 
the present analysis allows general functionals, such as the integral of 
the {\em natural energy} of the physical system.
Furthermore, this work confirms that the theory developed in 
Acquistapace {\em et~al.}~[Adv.~Differential Equations, 2005]---crucially 
utilized here---encompasses widely differing PDE problems, from thermoelastic 
systems to models of acoustic-structure and, now, fluid-structure 
interactions.
\vspace{2mm} 
}

\section{Introduction}
In this paper we consider the optimal control problem with quadratic functionals
for a fluid-structure interaction model.   
Of major concern is well-posedness of the Riccati equations arising in
the minimization problem, along with the feedback synthesis of the (unique) 
optimal control. 
The fluid-structure interaction is modeled by a system of coupled partial differential
equations (PDE) comprising a Stokes system (the fluid) and a three-dimensional system of dynamic elasticity (the solid). The coupling occurs at an interface 
separating two regions occupied, respectively, by the fluid and the solid.
It is worth mentioning at the outset that it is assumed that the motion of the 
solid is due to infinitesimal displacements. 
Accordingly, the fluid-solid interface is {\em stationary};
this and other modeling issues are discussed, e.g., in \cite{du-etal}.
The mathematical description of the PDE system, that is the boundary control problem
\eqref{e:stokeslame-0}, as well as further literature will be given in the next Section.
Our main goal is to establish the validity of a Riccati theory that would allow to 
control the structure, via boundary controls, acting as forces applied to the 
interface. 

It is well known that---even in the case of a single PDE---one of the main difficulties 
in a rigorous derivation of the feedback synthesis of the optimal control 
is the presence of {\em boundary} controls (or, more generally, {\em unbounded} control 
actions), combined with the lack of smoothing effects propagated by the dynamics
(see, e.g., \cite{bddm} and \cite{las-trig-books}).
In fact, while the linear-quadratic control problem with unbounded control operator 
has a complete solution in the case of PDE models whose free dynamics is governed by 
an {\em analytic} semigroup, this solution may be out of reach in the case of other 
kind of dynamics.
In particular, the case of purely hyperbolic PDE with boundary/point control 
is peculiarly different\footnote{
In the infinite time horizon case the so called {\em gain} (or feedback) operator 
$\sB^*P$ is intrinsically unbounded and the analysis of the algebraic Riccati equations 
is subtle; see \cite{flandoli-etal}, \cite{zwart-weiss}, and the subsequent improvements 
in \cite{trig-contemporary}, \cite{barbu-las-trig}.} 
from the parabolic case:
it would suffice to recall that in finite time horizon problems the Riccati operator (or optimal cost operator) $P(t)$ does {\em not} satisfy the differential Riccati equations, unless the observation operator possesses a suitable smoothing property.
 
On the other hand, certain interconnected PDE systems combining parabolic and
hyperbolic effects may give rise to an {\em abstract} control system $y'=Ay+By$ 
which yields a {\em singular estimate} for the operator $e^{At}B$, near $t=0$. 
This property---which is an intrinsic feature of control systems ruled by 
analytic semigroups---has been first identified in the analysis of an 
acoustic-structure interaction (where the overall semigroup was {\em not} 
analytic); see \cite{avalos-las-0}.
The essence of such estimates is the following: the parabolic component does induce a singular estimate (as a consequence of analyticity of the corresponding semigroup),
while hyperbolicity `transports' this estimate across the system through the coupling. 
Thus, if one can show that a singular estimate is valid for the entire system, 
then the theory in \cite{las-cbms,las-trento} ensures a feedback control law 
with {\em bounded} (in the state space) gain operator, along with well-posed Riccati equations. 
This theory has been successfully applied to diverse composite PDE models, including 
some thermoelastic systems, beside to various acoustic-structure interactions.
Several illustrations are contained in \cite{las-trento} and \cite{las-trig-se-2};
see also \cite{b-l-t}, \cite{bucci-las-thermo}, and the recent \cite{bucci-jee07}.

For the fluid-structure interaction under investigation, which comprises a parabolic 
and a hyperbolic PDE, it was shown in \cite{las-tuff-3} that a singular estimate 
(for the corresponding abstract evolution) is satisfied in the {\em finite energy space}, 
as long as the penalization in the quadratic functional does not involve the 
{\em mechanical energy} at a truly {\em energy level}.  
More precisely, the study in \cite{las-tuff-3} established specific singular 
estimates and hence well-posedness of the Riccati equations in the special case of
penalization of the mechanical variables below the energy level (say, 
{\em sub-critial} penalization), yet allowing full penalization of the fluid variable.

The situation becomes much more difficult when the mechanical variables are 
penalized at the {\em critical} level of the energy (see the functional 
\eqref{functional:integral-of-energy}). 
In fact, not only the regularity results of \cite{las-tuff-2} do not apply, but
the theory pertaining to control systems which yield singular estimates
(\cite{las-trig-se-2}, \cite{las-tuff-2}) is no longer valid. 
(Indeed, if it were so, the gain operator would be bounded on the state space,
while we will show that this is not the case; see Remark~\ref{r:unbounded-gain}.)
   
The present work addresses the issue of solvability of the optimal control problem 
with {\em general} quadratic functionals (i.e.~including {\it critical} penalization) 
for the PDE model \eqref{e:stokeslame-0}.
As we shall see, we provide (a positive) answer to the question remained open
in \cite[Remark~6.1]{las-tuff-3}.
This will follow in light of the theory introduced in \cite{abl-2}, which is shown 
to cover the present case in view of the set of trace regularity results established 
and collected in Theorem~\ref{t:main}. 
The theory developed in \cite{abl-2} is more effective in capturing the relevant 
properties of the dynamics, especially the ones which emanate from hyperbolicity.
These ultimately allow to define the gain operator as an {\it unbounded} operator on a suitable functional space. 
In this respect, the variational aspect of the minimization process 
is critical in order to justify the arguments leading to well-posedness of Riccati equations.     

Let us recall that the optimal control theory in \cite{abl-2}, while relaxing 
the `singular estimate requirement', postulates other regularity conditions of global 
nature. 
This makes it possible to obtain meaningful solutions to the differential 
Riccati equations, despite the gain operator is not bounded on the state space.
This, however, does not affect the synthesis, as the optimal solution still belongs
to the domain of the gain operator. 
Originally arisen in the study of boundary control problems for an established system of  thermoelasticity (\cite{abl-1}), so far this theory has been shown to apply as well
in the case of certain acoustic-structure interaction model including thermal effects 
(\cite{bucci-applicationes}).

\smallskip 
The paper is organized as follows.
In Section~\ref{section:2} we introduce the boundary control problem under
investigation, along with the statements of our main results, namely 
Theorem~\ref{t:main-0} and Theorem~\ref{t:main}. 
Moreover, we briefly record some necessary notation 
and the fundamental well-posedness result pertaining to the uncontrolled 
counterpart of the PDE system.
Section~\ref{section:3} is entirely devoted to the proof of Theorem~\ref{t:main},
which establishes the novel, distinct boundary regularity properties 
(of the solutions to the PDE system) which will ultimately result in
solvability of the optimization problem, i.e. Theorem~\ref{t:main-0}.
Section~\ref{section:4} contains the proof of Theorem~\ref{t:main-0},
based upon the application of the theory in \cite{abl-2}.
Finally, a short Appendix collects the statements of the regularity 
results pertaining to the elastic component of the system---recently obtained 
in \cite{barbu-etal-1} and \cite{las-tuff-3}---which are crucially 
utilized in the proof of Theorem~\ref{t:main}. 


\section{The PDE model, statement of main results} \label{section:2}

{\bf The PDE model. }
The PDE model under investigation describes the interaction of a (very slow) viscous, incompressible fluid, with an elastic body in a three dimensional bounded domain.
Although the introduction of such models dates back to \cite{lions},
their PDE analysis has increased significantly only in the past decade. 
A mathematical description of the composite PDE system is
given below.
By $\Omega_f$ and $\Omega_s$ we denote the open smooth domains 
occupied by the fluid and the solid, respectively.
Then $\Omega\subset \mathbb{R}^3$ denotes the entire solid-fluid region, 
that is $\Omega$ is the interior of $\overline{\Omega}_f\cup\overline{\Omega}_s$.
The boundary of $\Omega_s$ is the {\em interface} between the fluid and the solid,
and is denoted by $\Gamma_s=\partial\Omega_s$.
We finally denote by $\Gamma_f$ the outer boundary of $\Omega_f$, namely
$\Gamma_f=\partial\Omega_f\setminus \partial\Omega_s$.
It is assumed that the motion of the solid is entirely due to infinitesimal 
displacements, and hence that the interface $\Gamma_s$ is {\em fixed}.

The velocity field of the fluid is represented by a vector-valued
function $u$, which satisfies a 
Stokes system in $\Omega_f$;
the scalar function $p$ represents, as usual, the pressure.
In the solid region $\Omega_s$ the displacement $w$ satisfies
the equations of linear elasticity.
(The density and the kinematic viscosity 
which usually appear in the Navier-Stokes equation are set equal to one, just to 
simplify the notation).
The coupling takes place on the interface $\Gamma_s$. 
We recall from \cite{du-etal} that the interface condition 
$u=w_t$ on $\Gamma_s$ (in place of the usual no-slip boundary condition
$u=0$) accounts for the fact that although the displacement of the elastic 
body is small, its velocity is not (small, yet rapid oscillations).
Thus, the PDE system is given by 

\begin{equation} \label{e:stokeslame-0}
\left\{ \hspace{1mm}
\begin{split} 
& u_t-{\rm div}\,\epsilon(u) 
+\nabla p= 0 & &\textrm{in }\; Q_f:= \Omega_f\times (0,T)
\\
& {\rm div}\, u=0 & & \textrm{in }\; Q_f
\\
& w_{tt} - {\rm div}\,\sigma(w)=0 & &\textrm{in }\; Q_s:= \Omega_s\times (0,T)
\\
& u=0 & & \textrm{on }\; \Sigma_f:= \Gamma_f\times (0,T)
\\
& w_t=u & & \textrm{on }\; \Sigma_s:= \Gamma_s\times (0,T)
\\
& \sigma(w)\cdot \nu=\epsilon(u)\cdot\nu - p\nu-g & & 
\textrm{on }\; \Sigma_s
\\
& u(0,\cdot)=u_0 & & \textrm{in }\; \Omega_f
\\
& w(0,\cdot)=w_0\,, \quad w_t(0,\cdot)=w_1 & & \textrm{in }\; \Omega_f\,.
\end{split}
\right.
\end{equation}
In the above coupled PDE system, $\sigma$ and $\epsilon$ denote the elastic stress 
tensor and the strain tensor, respectively, that are
\begin{equation}
\sigma_{ij}(u) = \lambda \sum_{k=1}^3 \epsilon_{kk}(u)\delta_{ij}
+ 2\mu \epsilon_{ij}(u)\,,
\qquad
\epsilon_{ij}(u)= \frac12\Big(\frac{\partial u_i}{\partial x_j}
+\frac{\partial u_j}{\partial x_i}\Big)\,,
\end{equation}
where $\lambda, \mu$ are the Lam\'e constants and $\delta_{ij}$ is the Kronecker
symbol.


\smallskip
Since the present work is focused on the optimization problem, 
the subtle questions related to the modeling of fluid-structure interaction 
phenomena, as well as to the analysis of well-posedness of the corresponding 
coupled PDE systems, will not be discussed here.
Yet, well-posedness of the boundary value problem \eqref{e:stokeslame-0}, 
with $g\equiv0$ (that is the uncontrolled system \eqref{e:stokeslame-free} below),
is a prerequisite for the study of the associated optimal control problems.
Thus, although many authors have contributed to the PDE analysis of {\em nonlinear} 
fluid-structure interaction models (where the dynamics of the fluid is ruled by a 
Navier-Stokes equation), existence of finite energy weak solutions---even for the 
simpler Stokes-Lam\'e system \eqref{e:stokeslame-free}---has been an open question 
until \cite[Theorem~2.2]{barbu-etal-1}.
The reader is referred to \cite{barbu-etal-1,barbu-etal-2} for the analysis 
of well-posedness of the coupled PDE system \eqref{e:stokeslame-free};
in addition, \cite{barbu-etal-2} includes a very clear introduction to the  
(nonlinear) fluid-structure interaction problem, along with a technical comparison 
with the previous mathematical literature.
In Section~\ref{sub-sec-2-1} we shall recall, for the reader's convenience, 
the theory in \cite{barbu-etal-1} that is needed for our purposes.

We finally note that while the present study follows the variational approach 
of \cite{barbu-etal-1}, exploiting the novel boundary regularity results established 
therein, semigroup well-posedness and stability properties of the {\em linear} 
model have been investigated in \cite{avalos-trig-1}; see the survey paper 
\cite{avalos-trig-overview} and its references.
For the uniform stabilization problem, see \cite{avalos-trig-3}.

\medskip
\noindent
{\bf Further references.}
There is a large literature on coupled fluid-structure evolution problems.
Most works address the issue of developing models for specific physical 
problems and/or their numerical simulation.
Two main different scenarios arise from the applications: the case in which
the fluid is flowing in a tube with elastic walls, such as the blood through 
arteries, and the case where one or more elastic bodies are immersed in a 
fluid flow.
The PDE model under investigation pertains to a physical situation falling
under the latter category.

A very nice introduction to fluid-structure interaction problems is provided 
by \cite{du-etal}.
Recent treatises with focus on modeling and numerical analysis are 
\cite{quarteroni-formaggia} and \cite{moubachir-zol}.
An in-depth PDE analysis of well-posedness of these nonlinear models has indeed
appeared only recently.
Relevant contributions to this problem are given (without any claim of completeness) 
by \cite{sanmartin-etal}, \cite{feireisl}, the aforesaid \cite{du-etal}, 
\cite{coutand-shkoller}, \cite{boulakia}, \cite{barbu-etal-1,barbu-etal-2}
and, lastly, \cite{kukavica-etal}.
For more information on this subject, see the bibliography therein.

\subsection{Variational and semigroup formulation} \label{sub-sec-2-1}
Before giving the statement of our main results, let us preliminarly recall from 
\cite{barbu-etal-1} some basic notation, and the chief facts which pertain to the 
uncontrolled problem, that is system \eqref{e:stokeslame-0} with $g\equiv 0$.
Further technical results obtained in \cite{barbu-etal-1} and \cite{las-tuff-3}
will be needed in the proof of our main result; these will be recorded 
in an Appendix for convenience.

\noindent 
{\bf The uncontrolled model.}
Let us introduce the {\em free} system corresponding to \eqref{e:stokeslame-0},
namely
\begin{equation} \label{e:stokeslame-free}
\left\{ \hspace{1mm}
\begin{split} 
& u_t-{\rm div}\,\epsilon(u) 
+\nabla p= 0 & &\textrm{in }\; Q_f:= \Omega_f\times (0,T)
\\
& {\rm div}\, u=0 & & \textrm{in }\; Q_f
\\
& w_{tt} - {\rm div}\,\sigma(w)=0 & &\textrm{in }\; Q_s:= \Omega_s\times (0,T)
\\
& u=0 & & \textrm{on }\; \Sigma_f:= \Gamma_f\times (0,T)
\\
& w_t=u & & \textrm{on }\; \Sigma_s:= \Gamma_s\times (0,T)
\\
& \sigma(w)\cdot \nu=\epsilon(u)\cdot\nu - p\nu& & 
\textrm{on }\; \Sigma_s
\\
& u(0,\cdot)=u_0 & & \textrm{in }\; \Omega_f
\\
& w(0,\cdot)=w_0\,, \quad w_t(0,\cdot)=w_1 & & \textrm{in }\; \Omega_f\,.
\end{split}
\right.
\end{equation}
The energy space for the PDE problem \eqref{e:stokeslame-free} is 
\begin{equation*}
Y = H\times H^1(\Omega_s)\times L_2(\Omega_s)\,,
\end{equation*}
where $H$ is defined as follows:
\begin{equation*}
H := \big\{ u\in L_2(\Omega_f): \; {\rm div}\, u=0\,, \, u\cdot \nu|_{\Gamma_f}=0\big\}\,.
\end{equation*}
In addition, we denote by $V$ the space defined as follows:
\begin{equation*}
V := \big\{ v\in H^1(\Omega_f): \; {\rm div}\, u=0\,, \, u|_{\Gamma_f}=0\big\}\,;
\end{equation*} 
we shall use the following distinct notation for the various inner products
which will occurr throughout the paper:
\begin{equation*}
(u,v)_f := \int_{\Omega_f}uv \,{\rm d}\Omega_f\,,\quad 
(u,v)_s := \int_{\Omega_s}uv \,{\rm d}\Omega_s\,,\quad
\langle u,v\rangle := \int_{\Gamma_s} uv \,{\rm d}\Gamma_s\,.
\end{equation*}
The space $V$ is topologized with respect to the inner product given by 
\begin{equation*}
(u,v)_{1,f} := \int_{\Omega_f}\epsilon(u)\epsilon(v) d\Omega_f\,;
\end{equation*}
the corresponding (induced) norm $|\cdot|_{1,f}$ is equivalent to the 
usual $H^1(\Omega_f)$ norm, in view of Korn inequality and the Poincar\'e inequality.
\begin{remark}
{\rm 
The norm $\|\cdot\|_{H^r(D)}$ in the Sobolev space $H^r(D)$
will be shortly denoted by $|\cdot|_{r,D}$ throughout the paper. 
Note that all the Sobolev spaces $H^r$ related to $u$ and $w$
are actually $(H^r)^3$: the exponent is omitted just for the 
sake of simplicity.
}
\end{remark}
Let us recall from \cite{barbu-etal-1} the definition of {\em weak} solutions to 
the (uncontrolled) PDE system \eqref{e:stokeslame-free}.

\begin{definition}[Weak solution] 
\label{def:weak-sol}
Let $(u_0,w_0,w_1)\in H$ and $T>0$. We say that a triple 
$(u,w,w_t)\in C([0,T],H\times H^1(\Omega_s)\times L_2(\Omega_s))$
is a weak solution to the PDE system \eqref{e:stokeslame-0} 
if
\begin{itemize}
\item
$(u(\cdot,0),w(\cdot,0),w_t(\cdot,0))=(u_0,w_0,w_1)$,
\item
$u\in L_2(0,T;V)$,
\item
$\sigma(w)\cdot \nu \in L_2(0,T;H^{-1/2}(\Gamma_s))$,
$\frac{d}{dt}w|_{\Gamma_s}=u|_{\Gamma_s}\in L_2(0,T;H^{1/2}(\Gamma_s))$, and
\item
the following variational system holds a.e. in $t\in (0,T)$:
\begin{equation} \label{e:variational-system}
\begin{cases}
\frac{d}{dt}(u,\phi)_{f} + (\epsilon(u),\epsilon(\phi))_{f}
 - \langle\sigma(w)\cdot \nu + g,\phi\rangle=0
\\[1mm]
\frac{d}{dt}(w_{t},\psi)_{s} + (\sigma(w),\epsilon(\psi))_{s}
-\langle\sigma(w)\cdot \nu,\psi\rangle=0\,,
\end{cases}
\end{equation}
for all test functions $\phi\in V$ and $\psi\in H^1(\Omega_s)$.
\end{itemize}
\end{definition}

\begin{remark}
{\rm 
It is important to emphasize that the regularity properties of the normal stresses 
(see the third item of Definition~\ref{def:weak-sol}) do not follow from the interior regularity of the fluid-structure variables. 
It is an independent regularity result, showing the exceptional behavior of 
hyperbolic traces. 
This regularity property is necessary in order to justify the variational definition 
of weak solutions (see \eqref{e:variational-system}).
While there are other definitions of solutions to nonlinear PDE models of fluid-structure interactions which do not require additional regularity on the boundary (see, e.g.,
\cite{lions}, \cite{du-etal}, \cite{avalos-trig-1}), yet these definitions are not 
adequate to variationally decouple the (finite energy) weak solutions of the two 
equations.
On the other hand, this decoupling is crucially important in the present analysis, 
aimed at identifying the distinctive regularity properties of the overall dynamics,
that play a major role in the study of the associated optimal control problems.
Exploiting the distinct features (analyticity and hyperbolicity) of the decoupled 
dynamics makes it possible to establish the sharpest results for the coupled PDE 
system. 
(This fact was recently utilized in \cite{kukavica-etal}, as well.)  
Consequently, the issue of ``hidden'' regularity of the hyperbolic component is 
central to the problem studied and its solution. 
}
\end{remark}

Existence of weak (global) solutions of a nonlinear generalization of the PDE 
problem \eqref{e:stokeslame-free} has been established in \cite{barbu-etal-1}.

\begin{theorem}[Existence of weak solutions, \cite{barbu-etal-1}]
Given any initial datum $(u_0,w_0,w_1) \in Y$ and any $T>0$, there exists a weak 
solution $(u,w,w_t)$ to the system \eqref{e:stokeslame-free} such that 
\begin{equation*}
\nabla w\big|_{\Gamma_s}\in L_2(0,T;H^{-1/2}(\Gamma_s))\,,
\qquad 
\frac{{\rm d}}{{\rm d}t} w\Big|_{\Gamma_s}=
w_t\big|_{\Gamma_s}\in L_2(0,T;H^{1/2}(\Gamma_s))\,.
\end{equation*}
\end{theorem}


\bigskip
\noindent
{\bf The control system, semigroup formulation. }
Aiming to apply the optimal control theory pertaining to a general class of 
evolutions---in the present case, the one developed in \cite{abl-2}---it is convenient 
to recast the boundary value problem \eqref{e:stokeslame-0} as an abstract control
system in a Hilbert space. 
Accordingly, let us introduce the fluid dynamic operator $A: V\to V'$, defined by 
\begin{equation*}
(Au,\phi) = -(\epsilon(u),\epsilon(\phi))
\qquad \forall \,\phi\in V\,,
\end{equation*}
and the (Neumann) map $N:L_2(\Gamma_s)\to H$ defined as follows:
\begin{equation*}
Ng=h \Longleftrightarrow
(\epsilon(h),\epsilon(\phi))= \langle g,\phi\rangle \qquad \forall \phi\in V\,.
\end{equation*}
The chief properties of the operators $A$ and $N$ are recalled in the Appendix.

Then, if we set $y=(u,w,w_t)$, the boundary value problem \eqref{e:stokeslame-0} 
reduces to the linear control system 
\begin{equation} \label{e:control-system}
\begin{cases}
y'= \sA y + \sB g & {\rm in } \quad [\sD({\sA}^*)]'
\\
y(0)=y_0
\end{cases}
\end{equation}
where the (dynamic) operator $\sA:\sD(\sA)\subset Y\to Y$ is defined by
\begin{equation}\label{sA}
\begin{array}{l}
\sA=
\left(
\begin{array}{ccc}
A & AN \sigma(\,\cdot\,)\cdot\nu & 0\\
0 & 0 & I\\
0 & {\rm div}\,\sigma(\cdot) & 0
\end{array}
\right)\,, \quad 
\end{array}
\end{equation}
with domain
\begin{equation*}
\begin{array}{l}
\sD(\sA)= \big\{y=(u,w,z)\in H: \, u\in V, \; A(u+N\sigma(w)\cdot\nu)\in H\,,\;
 z\in H^1(\Omega_s)\,,
\\[2mm]
\qquad\qquad\qquad {\rm div}\, \sigma(w)\in L_2(\Omega_s)\,, \; 
z|_{\Gamma_s}=u|_{\Gamma_s}\big\}\,,
\end{array}
\end{equation*}
and the (control) operator $\sB:L_2(\Gamma_s)=U\to [\sD(\sA)]'$ reads as
\begin{equation}\label{sB}
\begin{array}{l}
\sB=
\left(
\begin{array}{c}
AN\\
0\\
0
\end{array}
\right)\,. 
\end{array}
\end{equation}
Given a quadratic functional
\begin{equation} \label{functional-general}
J(g) = \int_0^T\big(|\sR y(t)|_Z^2 + |g(t)|_U^2\big)\,dt\,,
\end{equation}
the optimal linear-quadratic (LQ) control problem 
is to minimize the functional \eqref{functional-general},
over all control functions $g\in L_2(0,T;U)$, with $y$ solution
to \eqref{e:control-system} corresponding to $g$.
As already pointed out in the Introduction, we aim to include in the present 
analysis non-smoothing observation operators $\sR$, such as the {\em identity}
operator; hence, $\sR$ is initially assumed to satisfy just $\sR\in \sL(Y,Z)$.
By doing so we admit natural quadratic functionals such as the following,
\begin{equation} \label{functional:integral-of-energy}
J(g) = \frac12\int_0^T\Big\{\,
|u(t)|_{0,\Omega_f}^2 +(\sigma(w(t)),\epsilon(w(t)))_{s}
+ |w_t(t)|_{0,\Omega_s}^2 + |g(t)|_{0,\Gamma_s}^2\,\Big\}\,dt
\end{equation}
which penalizes the full quadratic energy $E(t)$ of the system.
\begin{remark}
{\rm 
We already emphasized that the study performed in \cite{las-tuff-3} did not provide 
solvability of optimal control problems with general quadratic functionals: 
in particular, it did not cover the case of natural functionals such as  
\eqref{functional:integral-of-energy}.
On the other hand, the analysis carried out in \cite{las-tuff-3}---despite the
final constraint on the observation operator $\sR$---included the case of Bolza 
problems, where the penalization affects also the state at the {\em final time} 
$T<\infty$, namely, when the functional to be minimized is given by
\begin{equation} \label{bolza-functional}
J(g) = \int_0^T\big(|\sR y(t)|_Z^2 + |g(t)|_U^2\big)\,dt + (\sG y(T),y(T))_W\,.
\end{equation}
Note that the LQ-problem with Bolza-type quadratical functionals is not discussed 
here.
In fact, the LQ-problem with quadratic functionals of the form \eqref{bolza-functional} 
(with $\sG\neq 0$) for the class of control systems \eqref{e:control-system} described 
by the Assumptions~\ref{h:key}, has not been investigated yet.
}
\end{remark}

\subsection{Statement of the main results}
The main result of the present work is the proof of well-posedness of the 
(differential) Riccati equations 
corresponding to the optimal control problems associated with the 
fluid-structure model \eqref{e:stokeslame-0}, along with all the inherent assertions  
about solvability of the optimization problem; see Theorem~\ref{t:main-0}.
This variational result, however, critically relies on the novel trace regularity
results established specifically for the (uncontrolled) PDE system \eqref{e:stokeslame-free}
in Theorem~\ref{t:main}, which thus constitute the major technical contribution of the present work.
As we shall see, the proof of this 
set of regularity results is based on the interplay 
between the maximal parabolic regularity of the fluid component with the `hidden' 
regularity of the traces of the hyperbolic (solid) component. 
Indeed, the fact that the coupling is of hyperbolic/parabolic type will be critically 
utilized.

\subsubsection{The solution to the optimization problem}
With reference to the PDE model introduced in the previous section, let us consider the  optimal control problem \eqref{e:control-system}--\eqref{functional-general},
that is 

\begin{itemize}
\item[]
Minimize the functional $J(g)$ in \eqref{functional-general}, over all 
$g\in L_2(0,T; L_2(\Gamma_s))$, where $y(\cdot)=y(\cdot;y_0,g)$ solves the 
control system \eqref{e:control-system}.
\end{itemize}

Then we have the following.

\begin{theorem}\label{t:main-0}
Consider the optimal control problem \eqref{e:control-system}--\eqref{functional-general},
with $\sA$ and $\sB$ 
given by \eqref{sA} and \eqref{sB}, respectively.
If the observation operator satisfies  
\begin{equation} \label{h:observation}
{\sR}^*\sR\in \sL(\sD({\sA}^{\epsilon}),\sD({{\sA}^*}^{\epsilon}))
\end{equation}
for some $\epsilon\in (0,1/4)$, then the following assertions hold true.

\begin{enumerate}
\item
For any initial state $y_0\in Y$ there exists a unique optimal control 
$g^0(\cdot)\in L_2(0,T;L_2(\Gamma_s)$ such that 
\begin{equation*}
J(g^0)=\min_{g \in L_2(0,T;L_2(\Gamma_s)} J(g)\,.
\end{equation*}
The optimal pair $(g^0(\cdot),y^0(\cdot))$ has the following additional regularity: 
\begin{align}\label{state-reg}
& y^0(\cdot)= [u^{0}(\cdot), w^{0}(\cdot), w^{0}_{t}(\cdot)] 
\in C([0,T];H\times H^{1}(\Omega_{s}) \times L_2(\Omega_{s}))\,;
\nonumber \\[1mm]
& g^0(\cdot) \in \bigcap_{1\le p < \infty} L_p(0,T;L_2(\Gamma_s))\,.
\end{align}

\item
There exists a non-negative, selfadjoint operator (the Riccati operator) 
$P(t)\in \sL(Y)$, $t\in [0,T]$, defined explicitly in terms of the data, such that 
\begin{equation*}
J(g^{0})= (P(0)y_0,y_0)_Y\,;
\end{equation*}
more precisely, $P(\cdot)\in \sL(Y,C([0,T],Y))$.

\item
The gain operator $\sB^*P(\cdot)$ satisfies 
$\sB^*P(\cdot)\in \sL(\sD(\sA^{\epsilon}),C([0,T],L_2(\Gamma_s)))$;
moreover, one has (the feedback synthesis of the optimal control) a.e. in $[0,T]$:
\begin{equation}\label{e:feedback-formula}
g^0(t) = -\sB^* P(t)y^0(t)\,, \qquad \forall y_0\in Y\,.
\end{equation}

\item
The operator $P(t)$ is the {\em unique} solution of the Differential Riccati 
Equation satisfied for $0\le t< T$ and $x =(x_1,x_2,x_3), y = (y_1,y_2,y_3 )\in D(\sA)$,
\begin{equation} \label{e:dre}
\frac{d}{dt}(P(t) x,y)_Y + (\sA^*P(t)x,y)_Y + (P(t)\sA x,y)_Y 
+ (\sR x,\sR y)_Y= (\sB^* P(t)x,\sB^*P(t)y)\,,
\end{equation}
with 
\begin{equation*} 
\lim_{t\rightarrow T^-} (P(t)x,x) = 0 \qquad  \forall x \in Y\,.
\end{equation*}

\end{enumerate}

\end{theorem}

\begin{remark}
{\rm 
Since here the dynamics operator $\sA$ is the generator of a 
s.c.~contraction semigroup with $\sA^{-1}\in \sL(Y)$, then
the domains of fractional powers $\sD(\sA^{\epsilon})$ 
in \eqref{h:observation} may be computed as intermediate spaces 
between $\sD(\sA)$ and $Y$.
The same holds for $\sD({\sA^*}^{\epsilon})$.
(For a comprehensive list of cases where the identity 
$[\sD(\sA),Y]_{1-\theta}=\sD(\sA^{\theta})$ holds true, 
see, e.g., \cite[\S~0.2.1]{las-trig-books}).
Then, it is not difficult to show that in the present case
$\sD({\sA^*}^{\epsilon})\equiv \sD(\sA^{\epsilon})$,
provided $\epsilon$ is sufficiently small.
Therefore, assumption \eqref{h:observation} is satisfied, with a non-smoothing 
observation operator, such as $\sR=I$. 
This natural situation was indeed left as an open problem in \cite{las-tuff-3}.
}
\end{remark}

\begin{remark} \label{r:unbounded-gain}
{\rm 
Observe that the optimal pair does not display the typical regularity (in time) exhibited 
in the case of control systems whose underlying semigroup is analytic (or, more in general, when singular estimates are satisfied).
In particular, the optimal control is not continuous.
This is not surprising, in view of the influence of both hyperbolic and parabolic effects 
on the overall behavior of the solutions. 

Moreover, the gain operator $\sB^*P(t)$ is no longer bounded on the state space $Y$, 
but just {\em densely defined}.
However, this does not affect the final result, as the feedback formula 
holds for any initial state in the finite energy space. 
Thus, the observation operator $\sR$ need not have regularizing effects, 
and $\sR$ can be critical.  

The above observations also explain why the previous Riccati theories 
are intrinsincally unapplicable in the critical case, as they lead to
{\it bounded} gain operators, in contrast with the case under examination.
}
\end{remark}

\subsubsection{Trace estimates} 
The fundamental 
analytic tool which will enable us to show Theorem \ref{t:main-0} is a
complex of boundary regularity results pertaining to the fuid component 
of the PDE problem \eqref{e:stokeslame-free}.
These traces' regularity estimates of $u$ (and $u_t$) on the interface $\Gamma_s$
are the ``PDE counterpart'' of the abstract regularity properties of 
the (unbounded) operator $\sB^*e^{{\sA}^*t}$ needed to invoke the optimal
control theory of \cite{abl-2}. 
These regularity estimates are, however, also of independent interest. 

\begin{theorem}[Traces' regularity] \label{t:main}
Consider the uncontrolled Stokes-Lam\'e system, namely the PDE system 
\eqref{e:stokeslame-free}.
Let $y(t)=(u(t),w(t),w_t(t))$ be the solution corresponding to 
initial data $y_0=(u_0,w_0,w_1)$. 
Then the fluid component $u$ admits a decomposition $u(t) = u_1(t)+u_2(t)$,
and the following statements pertain to the regularity of the traces of 
$u_1$, $u_2$ and $u_t$ on $\Gamma_s$, respectively.
\begin{description}
\item{(i)}
The component $u_1$ satisfies a pointwise (in time) ``singular estimate'', 
namely there exists a positive constant $C_T$ such that
\begin{equation}\label{e:singular-estimate}
|u_1(t)|_{L_2(\Gamma_s)} \le 
\frac{C_T}{t^{1/4+\delta}}|y_0|_Y \qquad \forall y_0\in Y\,, \quad 
\forall t\in (0,T]
\end{equation}
(with arbitrarily small $\delta>0$).
\item{(ii)}
The component $u_2$ satisfies the following regularity:
\begin{description}
\item{(iia)}
if $y_0\in Y$, then $u_2|_{\Gamma_s}\in L_p(0,T;L_2(\Gamma_s))$
for all (finite) $p\ge 1$;

\medskip

\item{(iib)}
if $y_0\in \sD(\sA^{\epsilon})$, $\epsilon\in (0,\frac14)$, then
$u_2|_{\Gamma_s}\in C([0,T],L_2(\Gamma_s))$.
\end{description}

\item{(iii)} 
Let now $y_0\in \sD(\sA^{1-\theta})$, with $\theta\in (0,\frac14)$.
Then, the fluid component $u$ of corresponding solution satisfies, 
for some $q\in (1,2)$,
\begin{equation}
u_t|_{\Gamma_s} \in L_q(0,T;L_2(\Gamma_s))
\end{equation}
continuously with respect to $y_0$, that is there exists a constant $C_T$
such that 
\begin{equation}  \label{ineq:trace-estimate}
\|u_t\|_{L_q(0,T;L_2(\Gamma_s))} \le C_T\, \|y_0\|_{\sD(\sA^{1-\theta})}\,.
\end{equation}
The exponent $q$ will depend on $\theta$: more precisely, given $\theta\in (0,\frac14)$, 
one has 
\begin{equation} \label{e:range-for-q}
1<q<\frac{4}{3+4\theta}\,.
\end{equation}

\end{description}
\end{theorem}
The remainder of the paper is devoted to the proof of the two main results stated in 
Theorem~\ref{t:main-0} and Theorem~\ref{t:main}. 
Section~\ref{section:3} deals with the above boundary regularity results, which 
will be next utilized in Section~\ref{section:4} to establish Theorem~\ref{t:main-0}.


\section{Proof of the trace regularity results} \label{section:3}
This section is entirely devoted to the proof of our main contribution, 
that is Theorem~\ref{t:main}.

\smallskip
\noindent
{\em Proof of Theorem~\ref{t:main} }
Our starting point is the equation satisfied by $u(\cdot)$, namely 
$u_t=Au+AN \sigma(w)\cdot\nu$, whose solutions are given by 
\begin{equation} \label{e:mild-1}
u(t) = e^{At}u_0+ \int_0^t e^{A(t-s)} AN \sigma(w)(s,\cdot)\nu \,ds\,;
\end{equation}
the above expression yields the natural splitting $u(t) =u_1(t)+u_2(t)$,
with 
\begin{equation}\label{e:terms-decomposition}
u_1(t):= e^{At}u_0\,, \qquad 
u_2(t):=\int_0^t e^{A(t-s)} AN \sigma(w)(s,\cdot)\nu \,ds\,.
\end{equation}
In view of $N^*Au=-u|_{\Gamma_s}$ (see Lemma~\ref{l:neumann-trace}), 
the corresponding traces on $\Gamma_s$ read as
\begin{align}
u_1|_{\Gamma_s} &=-N^*Au_1(t)= -N^*A e^{At}u_0\,, &
\nonumber\\
u_2|_{\Gamma_s} &=-N^*Au_2(t)=
-N^*A\int_0^t e^{A(t-s)} AN \sigma(w)(s,\cdot)\nu \,ds\,, &
\label{e:traces-expression}
\end{align}
respectively.

\noindent
{\em (i)}
The {\em singular} estimate in \eqref{e:singular-estimate}
follows as an immediate consequence of the well known estimates
pertaining to analytic semigroups:
\begin{equation}
\begin{split}
\big|u_1(t)|_{\Gamma_s}\big|_U= |N^*A e^{At}u_0|
=|N^*A^{3/4-\delta}A^{1/4+\delta} e^{At}u_0| 
\\[2mm]
\le \|N^*A^{3/4-\delta}\|\, |A^{1/4+\delta} e^{At}u_0|_Y
\sim C_Tt^{-1/4-\delta}|u_0|\,.
\end{split}
\end{equation}
This shows the validity of assertion $(i)$.

\smallskip
\noindent
{\em (ii)}
Let initially $y_0=(u_0,w_0,w_1)\in Y$.
In view of \eqref{e:traces-expression}, it is clear that the sharp regularity theory 
pertaining to the the wave component will play a central role in the study of the 
regularity of the traces of $u_2(t)$ on $\Gamma_s$. 
More precisely, we shall utilize the recent trace results obtained in 
\cite[Theorem~3.3]{barbu-etal-1} and refined in \cite[Lemma~5.2]{las-tuff-3}; 
see Lemma~\ref{l:wave-regularity} in the Appendix.
Accordingly, following the decomposition of $\sigma(w)\cdot \nu$ established
in Lemma~\ref{l:wave-regularity}, it is convenient to introduce a further splitting, namely 
\begin{equation}
u_2(t)=
\underbrace{\int_0^t e^{A(t-s)} AN \sigma(w_1)(s,\cdot)\nu \,ds}_{u_{21}(t)}+\underbrace{\int_0^t e^{A(t-s)} AN \sigma(w_2)(s,\cdot)\nu \,ds}_{u_{22}(t)}\,.
\end{equation}
Thus, one has first
\begin{equation*}
\begin{split}
N^*A u_{21}(t) 
&= N^*A \int_0^t e^{A(t-s)} AN \sigma(w_1)(s)\cdot\nu \,ds
\\
&= [N^*A^{3/4-\epsilon}] A^{1/4+\epsilon+1/2}
\int_0^t e^{A(t-s)} \underbrace{A^{1/2}N \sigma(w_1)(s)\cdot\nu}_{f(s)}\,ds
\end{split}
\end{equation*}
where $f\in C([0,T],L_2(\Omega_s))$ in view of Lemma~\ref{l:wave-regularity}
and Lemma~\ref{l:neumann-trace}.
Consequently,
\begin{equation} \label{e:continuity}
\int_0^t e^{A(t-s)} f(s)\,ds\in C([0,T],D(A^{1-\sigma}))\,,
\end{equation}
with arbitrarily small $\sigma>0$; see, e.g., \cite[Proposition~0.1, p.~4]{las-trig-books}.
Therefore $N^*A u_{21}\in C([0,T],L_2(\Gamma_s))$, and 
{\em a fortiori} we obtain 
\begin{equation}\label{e:leadsto(iia)-1}
N^*A u_{21}\in L_p(0,T;U) \qquad \forall p\ge 1\,.
\end{equation}

As for the second summand $N^*A u_{22}(t)$, because of the different regularity 
of $\sigma(w_2)(s)\cdot\nu$ we rewrite in a different fashion:
\begin{equation*}
N^*A u_{22}(t) = 
\big[N^*A^{3/4-\epsilon}\big]\, A^{1/2+2\epsilon}
\int_0^t e^{A(t-s)} \underbrace{[A^{3/4-\epsilon}N
\sigma(w_2)(s)\cdot\nu]}_{\varphi(s)} \,ds\,.
\end{equation*}
Notice now that the above integral is the convolution 
\begin{equation*}
\int_0^t K(t-s) \varphi(s)\,ds\,,
\end{equation*}
with $\varphi\in L_2(0,T;U)$ and the kernel $K$ such that 
$||K(s)|| \sim \frac1{s^{1/2+2\epsilon}}$, where $\epsilon$ can be taken
arbitrarily small.
Hence $K\in L_{2-\sigma}(0,T;U)$ for arbitrarily small $\sigma>0$.
Thus, the Young inequality yields 
\begin{equation}\label{e:leadsto(iia)-2}
N^*A u_{22}\in L_p(0,T;U) \qquad \forall p\ge 1\,.
\end{equation}
Thus, \eqref{e:leadsto(iia)-2} combined with \eqref{e:leadsto(iia)-1}
shows the validity of assertion $(iia)$.

\smallskip
Let now $y_0\in D({\cal{A}}^\epsilon)$, $\epsilon>0$.
In this case we have 
\begin{equation*}
D({\cal{A}}^\epsilon)\subset H\times H^{1+\epsilon}(\Omega_2)
\times H^{\epsilon}(\Omega_2)
\end{equation*}
and by Lemma~\ref{l:key} it follows
$u|_{\Gamma_s}\in H^\epsilon(\Sigma_s)$,
provided that $\epsilon<\frac14$.
This enables us to apply the second part of Lemma~\ref{l:wave-regularity},
which gives
\begin{equation}\label{e:additional}
\sigma(w_1)\cdot \nu\in C([0,T],H^{-1/2}(\Gamma_s))\,, \quad 
\sigma(w_2)\cdot \nu\in H^{\epsilon}(\Sigma_s)\,.
\end{equation}
Now, the analysis of $N^*Au_{21}$ follows closely the one in item $(iia)$, yielding the conclusion $N^*Au_{21}\in C([0,T],L_2(\Gamma_s))$
(this is justified by the membership \eqref{e:continuity}).
Instead, on the basis of the novel regularity of $\sigma(w_2)\cdot \nu$ 
in \eqref{e:additional}, from parabolic theory it follows that
\begin{equation*}
u_{22}\in H^{\epsilon+3/2,\epsilon/2+3/4}(Q_f)\,, 
\end{equation*}
so that 
\begin{equation*}
N^*Au_{22}\in H^{\epsilon+1,\epsilon/2+1/2}(\Sigma_s)
\subset H^{1/2+\epsilon/2}(0,T;L_2(\Gamma_s))\subset C([0,T],L_2(\Gamma_s))\,.
\end{equation*}
%
As both $N^*Au_{21}$ and $N^*Au_{22}$ belong to 
$C([0,T],L_2(\Gamma_s))$, then
$N^*Au_2\in C([0,T],L_2(\Gamma_s))$ and $(iib)$ is proved.

\smallskip
\noindent
{\em (iii)}
In this last step we aim to ascertain the regularity of the boundary traces 
of the time derivative $u_t$ on $\Gamma_s$.
We return to the mild solution \eqref{e:mild-1} and compute
\begin{equation}\label{e:decomp-1}
u_t(t) = \underbrace{Ae^{At}u_0}_{v_1(t)}
+ \underbrace{\overbrace{A\int_0^t e^{A(t-s)} AN \sigma(w)(s,\cdot)\nu\,ds}^{v_{2a}}
+ \overbrace{AN \sigma(w)(t,\cdot)\nu}^{v_{2b}}}_{v_2(t)}
\end{equation}
which can also be rewritten as
\begin{equation}\label{e:decomp-2}
u_t(t) =  \underbrace{Ae^{At}u_0}_{v_1(t)}
+ \underbrace{\overbrace{\int_0^t e^{A(t-s)} AN \sigma(w_s)(s,\cdot)\nu \,ds}^{v_{21}(t)}
+ \overbrace{Ae^{At}N \sigma(w)(0,\cdot)\nu}^{v_{22}(t)}}_{v_2(t)}\,.
\end{equation}
%
The plan we aim to carry out is to discuss first the regularity of 
the function $v_2:=\frac{\partial u_2}{\partial t}$ when $y_0\in D(\sA)$
by using its expression in \eqref{e:decomp-2}.
Next, when $y_0\in Y$, we would rather utilize \eqref{e:decomp-1}, and 
then use interpolation arguments to establish the regularity
corresponding to initial data in $D(\sA^{1-\epsilon})$.
Only subsequently we shall derive the trace regularity of $v_2$ by applying 
the operator $-N^*A$.

\smallskip
When $y_0\in D(\sA)$, by standard semigroup arguments we known that 
$\sigma(w_t)\cdot \nu$ exhibits the same regularity as that of 
$\sigma(w)\cdot\nu$ when $y_0\in Y$, i.e.
(invoking once again Lemma~\ref{l:wave-regularity})
\begin{equation*}
\sigma(w_t)\cdot \nu = \sigma_1 + \sigma_2 
\in C([0,T],H^{-1/2}(\Gamma_s))\,{\oplus}\,L_2(0,T;L_2(\Gamma_s))\,.
\end{equation*}
To pinpoint the regularity of $v_{21}$, we now utilize the above splitting and 
follow the analysis carried out in the proof of $(ii)$.
More precisely, combining elliptic regularity (in particular, Lemma~\ref{l:neumann-trace}), with the analyticity of the semigroup $e^{At}$, along with the (singular) 
estimates pertaining to $A^\alpha e^{At}$, we first obtain, for any $t$ and any 
$\delta<1/2$,  
\begin{multline}
\big|A^{1/2 -\delta}\int_0^t e^{A (t-s) } A N \sigma_1(s)\, ds\big|=
\big|A^{1 -\delta}\int_0^t e^{A (t-s) } A^{1/2} N \sigma_1(s)\, ds\big| 
\\[1mm]
\le C \int_0^t \frac{1}{(t-s)^{1-\delta}} ds \,\|\sigma_1\|_{C([0,T],H^{-1/2}(\Gamma_s))}
\le C \|\sigma_1\|_{C([0,T],H^{-1/2}(\Gamma_s))}\,.
\end{multline}
\\
As for the latter term, we apply as well $A^{1/2 -\delta}$ and rewrite as follows:
\begin{equation}
A^{1/2 -\delta} \int_0^t e^{A(t-s)}A N \sigma_2(s)\, ds 
=
\int_0^t [A^{3/4 -\delta/2} e^{A(t-s)}]\,[A^{3/4-\delta/2} N]\,\sigma_2(s)\, ds\,,
\end{equation} 
where it is clear now that the integral is the convolution of 
$L_{4/3}$ and $L_2$ (in time) functions, respectively.
On the strength of the Young's inequality, we get $L_4$-regularity in time, 
so that 
\begin{equation*}
v_{21}\in C([0,T],D(A^{1/2-\delta}))\,{\oplus}\,L_{4}(0,T;D(A^{1/2-\delta}))\,,
\qquad 0<\delta<\frac12\,.
\end{equation*}
This implies the membership
\begin{equation}\label{e:m1}
v_{21}\in L_4(0,T;D(A^{1/2-\delta}))\,, \qquad 0<\delta<\frac12\,.
\end{equation}
On the other hand, still with $y_0\in D(\sA)$, one has just 
$\sigma(w)\cdot\nu\in C([0,T],H^{-1/2}(\Gamma_s))$ 
which suggests us to rewrite $v_{22}$ as follows:
\begin{equation*}
v_{22}(t)= A^{1/2}e^{At}\,\big(A^{1/2}N \sigma(w)(0,\cdot)\nu\big)\,;
\end{equation*}
then, again in view of Lemma~\ref{l:neumann-trace} and of the usual 
singular estimates pertaining to analytic semigroups, it follows
\begin{equation}\label{e:m2}
v_{22}\in L_q(0,T;D(A^{1/2-\delta}))\,, 
\end{equation}
provided that $q(1-\delta)<1$. 
Therefore, \eqref{e:m1} combined with \eqref{e:m2} yields, for any $0<\delta_1<\frac12$
\begin{equation}\label{e:m}
y_0\in D(\sA) \Longrightarrow v_2\in L_{q_1}(0,T;D(A^{1/2-\delta_1}))
\equiv L_{q_1}(0,T;H^{1-2\delta_1}(\Omega_f))\,,
\end{equation}
where $q_1\in (1,2)$ depends on $\delta_1$; more precisely, 
\begin{equation}\label{e:sobolev-exponent_1}
q_1< \frac1{1-\delta_1}\,. 
\end{equation}

\smallskip
Let now $y_0\in Y$. 
In this case we use the decomposition \eqref{e:decomp-1}, and begin with the 
analysis of $v_{2a}$.
Setting $w=w_1+w_2$ (according with Lemma~\ref{l:wave-regularity}), one has
\begin{align*}
& A\int_0^t e^{A(t-s)} AN \sigma(w_1)(s,\cdot)\nu\,ds=
\\
& \quad A^{1/2+\epsilon_1}A^{1-\epsilon_1}\int_0^t e^{A(t-s)} A^{1/2}N \sigma(w_1)(s,\cdot)\nu\,ds \in C([0,T],[D(A^{1/2+\epsilon_1})]')\,,
\end{align*}
while
\begin{align*}
& A\int_0^t e^{A(t-s)} AN \sigma(w_2)(s,\cdot)\nu\,ds=
\\
& \quad A^{1/4+\epsilon_2}A\int_0^t e^{A(t-s)} A^{3/4-\epsilon_2}N \sigma(w_2)(s,\cdot)\nu\,ds
\in L_2(0,T;[D(A^{1/4+\epsilon_2})]') 
\end{align*}
where both $\epsilon_1$ and $\epsilon_2$ can be taken arbitrarily small.
As a result, 
\begin{equation}\label{e:regularity-v_2a}
v_{2a}\in L_2(0,T;[D(A^{1/2+\epsilon})]')\,, \qquad 0<\epsilon<\frac12\,.
\end{equation}
As for the term $v_{2b}$, readily 
\[
AN\sigma(w_1)(t,\cdot)\nu= A^{1/2} A^{1/2}N\sigma(w_1)(t,\cdot)\nu
\in C([0,T],[D(A^{1/2})]')
\]
while 
\[
AN\sigma(w_2)(t,\cdot)\nu= A^{1/4+\epsilon} [A^{3/4-\epsilon}N]\sigma(w_2)(t,\cdot)\nu
\in L_2(0,T;[D(A^{1/4+\epsilon})]')\,, 
\]
and since $\epsilon$ can be taken arbitrarily small, we deduce as well
\begin{equation}\label{e:regularity-v_2b}
v_{2b}=AN\sigma(w)(t,\cdot)\nu\in L_2(0,T;[D(A^{1/2})]')\,.
\end{equation}
On the basis of \eqref{e:regularity-v_2a} and \eqref{e:regularity-v_2b},
we obtain
\begin{equation}\label{e:m-0}
y_0\in Y \Longrightarrow v_2=v_{2a}+v_{2b}\in L_2(0,T;[D(A^{1/2+\delta_2})]')
\equiv L_2(0,T;[H^{1+2\delta_2}(\Omega_f)]')\,,
\end{equation}
if $0<\delta_2<\frac14$.

\smallskip
Thus, \eqref{e:m-0}, combined with \eqref{e:m}, gives by interpolation
\begin{equation}\label{e:interpolation}
y_0\in D(\sA^{1-\theta})\Longrightarrow  v_2\in L_{q_1}(0,T;W)\,, 
\end{equation}
where $q_1$ is as in \eqref{e:sobolev-exponent_1} and $W$ is the interpolation space 
\[
W= (H^{1-2\delta_1}(\Omega_f),[H^{1+2\delta_2}(\Omega_f)]')_{\theta}
\equiv H^{s}(\Omega_f)\,, 
\]
if
\begin{equation*}
s=(1-2\delta_1)(1-\theta)-\theta(1+2\delta_2)=
1-2\delta_1-2\theta(1+\delta_2-\delta_1)\ge 0\,;
\end{equation*}
see \cite[Theorem~12.5]{lions-magenes}.
Notice that by taking, for instance, $\delta_1=\delta_2=:\delta$, 
one has $s\ge 1/2$ provided that
\begin{equation}\label{e:theta-constraint}
\theta+\delta \le \frac14\,.
\end{equation}
In this case $v_2\in H^s(\Omega_f)$ with $s\ge 1/2$ and hence its trace on $\Gamma_s$ 
is well defined. 
Notice that, in view of the constraint \eqref{e:theta-constraint},
we need to require $0<\theta<\frac14$.
Consequently, given any $\theta$ such that $0<\theta<\frac14$, 
choosing, e.g., $\delta= 1/4 -\theta$ in view of \eqref{e:theta-constraint}, 
from \eqref{e:interpolation} it follows
\begin{equation}\label{e:trace_1}
y_0\in D(\sA^{1-\theta})) \Longrightarrow
N^*Av_2\in L_{q_1}(0,T;L_2(\Gamma_s)) \qquad \forall q_1< \frac{4}{3+4\theta}\,.
\end{equation}

It remains to establish the regularity of the first summand
$N^*Av_1(t)= N^*Ae^{At}Au_0$ when $y_0\in D(\sA^{1-\theta})$.
In this case $u_0\in (H^1(\Omega_f),L_2(\Omega_f))_\theta=H^{1-\theta}(\Omega_f)$, 
and from
\begin{equation*}
N^*Av_1(t):=N^*Ae^{At}Au_0=[N^*A^{3/4-\epsilon}] A^{1/4+\epsilon+1/2+\theta/2}e^{At}
A^{(1-\theta)/2}u_0\,,
\end{equation*}
it immediately follows
\begin{equation}\label{e:trace_2}
y_0\in D(\sA^{1-\theta}) \Longrightarrow
N^*Av_1\in L_{q_2}(0,T;L_2(\Gamma_s)) \qquad \forall q_2< \frac{4}{3+2\theta+4\epsilon}\,.
\end{equation}
Notice that in the above membership the Sobolev exponent $q_2$ belongs to 
$(1,2)$, as well.
In conclusion, since $\epsilon$ in \eqref{e:trace_2} can be taken arbitrarily small,
the regularity $L_{q_1}(0,T;L_2(\Gamma_s))$ combined with $L_{q_2}(0,T;L_2(\Gamma_s))$ 
(in \eqref{e:trace_1} and \eqref{e:trace_2}, respectively) imply the membership 
\begin{equation}\label{e:trace}
y_0\in D(\sA^{1-\theta}) \Longrightarrow
u_t|_{\Gamma_s}=:v|_{\Gamma_s}
\in L_{q}(0,T;L_2(\Gamma_s)) 
\qquad \forall q< \frac{4}{3+4\theta}\,,
\end{equation}
which concludes the proof.


\section{Proof of Theorem~\ref{t:main-0}} \label{section:4}
The conclusions stated in Theorem~\ref{t:main-0} will follow from 
\cite[Theorem~2.3]{abl-1}, once we verify the standing assumptions,
which are recorded below for the reader's convenience.
\begin{assumptions}\label{h:key}
For each $t\in [0,T]$, the operator $B^*e^{tA^*}$ can be represented as 
\begin{equation} \label{keyass} 
{\sB}^*e^{{\sA}^*t}y_0 = F(t)y_0 + G(t)y_0, \qquad t\ge 0, \quad y_0\in \sD({\sA}^*)\,,
\end{equation}
where $F(t):Y\to U$ and $G(t):\sD({\sA}^*)\to U$, $t>0$, are bounded linear
operators satisfying the following assumptions:

\begin{enumerate}

\item
there is $\gamma\in \, (\frac{1}{2},1)$ such that 
$\|F(t)\|_{\sL(Y,U)} \le C\, t^{-\gamma}$ for all $t\in (0,T]$;

\item 
the operator $G(\cdot)$ belongs to $\sL(Y,L^p(0,T;U))$ for all $p\in [1,\infty)$, with 
\begin{equation} \label{keyasGa}
\|G(\cdot)\|_{\sL(Y,L^p(0,T;U))} \le c_p < \infty \qquad \forall p\in [1,\infty)\,;
\end{equation}

\item 
there is $\varepsilon>0$ such that:

\begin{enumerate}

\item
the operator $G(\cdot){{\sA}^*}^{-\varepsilon}$ belongs to
$\sL(Y,C([0,T],U))$, with
\begin{equation} \label{keyasGb}
\|{\sA}^{-\varepsilon}G(t)^*\|_{\sL(U,Y)} \le c <\infty \qquad \forall t\in [0,T]\,;
\end{equation}

\item
the operator ${\sR}^*\sR$ belongs to
$\sL(\sD({\sA}^{\varepsilon}),\sD({{\sA}^*}^{\varepsilon}))$;

\item 
there is $q\in\, (1,2)$ (depending, in general, on $\varepsilon$) such that the operator 
${\sB}^*e^{{\sA}^*\cdot}{\sR}^*\sR{\sA}^{\varepsilon}$ has an 
extension belonging to $\sL(Y,L^q(0,T;U))$. 
\end{enumerate}
\end{enumerate}

\end{assumptions}

\begin{remark}
{\rm 
Note that the set of requirements in Assumptions~\ref{h:key} involves the regularity (in time) of the operator ${\sB}^*e^{{\sA}^*t}$, both locally at the origin (with singularity controlled by $\gamma$), and globally (in $L_p$ sense). 
}
\end{remark}

We now prove that the regularity results established in Theorem~\ref{t:main} 
imply all the requirements in Assumptions~\ref{h:key}, with suitable 
values of $\gamma$, $\epsilon$ and $q$. 

\smallskip
\noindent
{\em Verification of Assumptions~\ref{h:key}. } 
Following \cite[Proof of Theorem~5.1]{las-tuff-3}, it is not difficult
to verify that given any initial state $y_0=(u_0,w_0,w_1)\in \sD(\sA^*)$, one has 
$\sB^* e^{\sA^*t}y_0 = N^*A \hat{u}(t)=- \hat{u}(t)|_{\Gamma_s}$, 
where $\hat{u}(t)$ is the first component of the solution 
$\hat{y}:=(\hat{u},\hat{w},\hat{w_t})$ to the (homogeneous) {\em adjoint} system 
\begin{equation*}
\begin{cases}
\hat{y}'(t)=\sA^*\hat{y}(t) \\
\hat{y}(0)=y_0\,.
\end{cases}
\end{equation*}
The abstract adjoint system yields a system of coupled PDE
which is essentially the same as the original boundary value problem 
\eqref{e:stokeslame-free}, except for few changes of sign in the equations.
As a consequence, the regularity results established by Lemma~\ref{l:wave-regularity}
and Lemma~\ref{l:key} readily produce analogues, by replacing $\sD(\sA)$ 
and $y=(u,w,w_t)$ by $\sD(\sA^*)$ and $\hat{y}=(\hat{u},\hat{w},\hat{w}_t)$, 
respectively.
Similarly, the fluid component $\hat{u}$ of the solution $\hat{y}$ to 
the dual PDE system 
satisfies---{\em mutatis mutandis}---the regularity properties in Theorem~\ref{t:main}.

\noindent
1. In light of the decomposition of $\hat{u}$ found in Theorem \ref{t:main}, 
let us set 
\begin{equation*}
F(t) y_0 :=\hat{u}_1(t)\big|_{\Gamma_s}\,, \qquad 
G(t) y_0 := \hat{u}_2(t)\big|_{\Gamma_s}\,.  
\end{equation*}
Then, the first statement in Theorem~\ref{t:main}, along with the estimate 
\eqref{e:singular-estimate}, provides us with the soughtafter singular 
estimate, with (optimal) exponent $\gamma = 1/4 + \delta$, 
and the first of Assumptions~\ref{h:key} is satisfied.
\\
2. Assertion {\em (iia)} in Theorem~\ref{t:main} is nothing but the regularity 
condition 2. of the Assumptions~\ref{h:key}, valid for all $p\in [1,\infty)$. 
\\
3. Condition {\em (iib)} of Theorem~\ref{t:main} translates into 
$G(t){\sA^*}^{-\epsilon} y_0 \in C([0,T],L_2(\Gamma_s))$,  
which in turn gives the assertion (3a) of the Assumptions~\ref{h:key}, 
with no constraints on $\epsilon$. 
It remains to verify the tricky assertion (3c) of Assumptions~\ref{h:key}. 
This will be implied by condition {\em (iii)} in Theorem~\ref{t:main}.
%
We first claim that the estimate \eqref{ineq:trace-estimate}
in {\em (iii)} of Theorem~\ref{t:main} yields, for any $\theta \in (0,1/4)$, 
the following one:
\begin{equation}\label{e:equivalent}
\|\sB^* e^{\sA^* t}{\sA^*}^\theta y_0\|_{L^q(0,T;L_2(\Gamma_s))} 
\le C |y_0|_Y\,, \quad y_0\in \sD({A^*}^{\theta})\,, 
\quad 1<q<\frac{4}{3 + 4 \theta}\,.
\end{equation}
This is easily seen if one just observes that if $y_0\in \sD({\sA^*}^\theta)$ one has
\begin{equation*}
\sB^*e^{\sA^*t}{\sA^*}^\theta y_0=\sB^*e^{\sA^*t}\sA^* \,({\sA^*}^{\theta-1}y_0)
=\sB^* \frac{d}{dt}\,e^{\sA^*t}z_0= \hat{u}_t|_{\Gamma_s} \,, 
\end{equation*}
where now $z_0:={\sA^*}^{\theta-1}y_0\in \sD({\sA^*}^{1-\theta})$.
Then, in view of the assumption (3b) on the observation operator $\sR$,
one concludes
\begin{align*}
\|\sB^* e^{\sA^* t} \sR\sR^* {\sA}^{\theta} y_0\|_{L^q(0,T; L_2(\Gamma_s))}
&=
\|\sB^* e^{\sA^* t} {\sA^*}^\theta {\sA^*}^{-\theta}\sR\sR^* 
{\sA}^{\theta} y_0\|_{L^q(0,T; L_2(\Gamma_s))}
\\[1mm]
& \le C |y_0|_Y\,, \quad  y_0\in \sD(\sA^{\theta})\,, 
\end{align*} 
i.e. condition (3c) is satisfied with $\epsilon = \theta$, $0<\theta<1/4$,
for any $q$ such that $1<q<4/(3 + 4 \theta)$.
This completes the proof of Theorem~\ref{t:main-0}.
\qed

\appendix

\section{Appendix}
For completeness' sake and for the reader's convenience we record 
some results shown in \cite{barbu-etal-1} and in \cite{las-tuff-3},
which are used frequently or critically in the proof of our main result.


\begin{lemma}[\cite{las-tuff-3}] 
\label{l:neumann-trace}
The Green (Neumann) map $N:L_2(\Gamma_s)\to H$ satisfies the following 
regularity results.
\\
$(i)$ One has $N^*Au = -u|_{\Gamma_s}$, $u\in V$, where the adjoint 
is computed with respect to the $L_2$-topology.
\\
$(ii)$ $N\in \sL(L_2(\Gamma_s),\sD(A^{3/4-\delta}))\cap
\sL(H^{-1/2}(\Gamma_s),\sD(A^{1/2}))$ for any $\delta$, $0<\delta<\frac34$. 
\end{lemma}


\begin{lemma}[\cite{barbu-etal-1,las-tuff-3}] 
\label{l:wave-regularity}
Let $(w_0,w_1)\in H^{\alpha+1}(\Omega_s)\times H^{\alpha}(\Omega_s)$,
with $0\le \alpha \le \frac14$, and let $f\in L_2(0,T;H^{1/2}(\Gamma_s))$.
Then, the solution of the initial/boundary value problem 
\begin{equation}
\begin{cases}
w_{tt}- {\rm div}\,\sigma(w) =0 & \textrm{in }\; Q_s
\\
\frac{d}{dt} w|_{\Gamma_s}=f & \textrm{on }\; \Sigma_s
\\
w(0,\cdot)=w_0\,,\; w_t(0,\cdot)=w_1\quad  & \textrm{in }\; \Omega_s
\end{cases}
\end{equation}
can be decomposed as $w=w_1+w_2$, where 
$\sigma(w_1)\cdot \nu\in C([0,T],H^{-1/2}(\Gamma_s))$, while
$\sigma(w_2)\cdot \nu\in L_2(0,T;L_2(\Gamma_s))$.
If, in addition, $f\in H^\alpha(\Sigma_s)$, then
$\sigma(w_2)\cdot \nu\in H^\alpha(\Sigma_s)$.
Moreover, the following estimates hold true.
\begin{align*}
\|\sigma(w_1)\cdot \nu\|_{C([0,T],H^{-1/2}(\Gamma_s))}^2
&\le C_1 \big(|w_0|_{1,\Omega_s}^2+|w_1|_{0,\Omega_s}^2
+|f|_{L_2(0,T;H^{1/2}(\Gamma_s))}\big)
\\[2mm]
\|\sigma(w_2)\cdot \nu\|_{H^{\alpha}(\Sigma_s))}^2
& \le C_2 \big(|w_0|_{1+\alpha,\Omega_s}^2+|w_1|_{\alpha,\Omega_s}^2
+|f|_{H^{\alpha}(\Sigma_s))}\big)
\end{align*}
\end{lemma}



\begin{lemma}[\cite{las-tuff-3}]
\label{l:key}
Consider the uncontrolled counterpart of the PDE problem~\eqref{e:stokeslame-0},
that is \eqref{e:stokeslame-0} with $g\equiv 0$.
Let initial data satisfy 
$(u_0,w_0,w_1)\in L_2(\Omega_f)\times H^{1+\alpha}(\Omega_s)
\times H^{\alpha}(\Omega_s)$, $0\le \alpha \le \frac14$.
Then, for any $T < \infty $ we have 
$u|_{\Gamma_s}\in H^{\alpha}(\Sigma_s)$ and the following estimate holds true:
\begin{equation}\label{e:improved}
|u|_{H^{\alpha}(\Sigma_s)}^2\le C 
\big(|u_0|_{0,\Omega_f}^2+|w_0|_{1+\alpha,\Omega_s}^2
+|w_1|_{\alpha,\Omega_s}^2\big)\,.
\end{equation}
\end{lemma}

\begin{remark} 
{\rm 
It it important to emphasize that the original proof of Lemma~\ref{l:key} 
(given in \cite[Lemma~5.3]{las-tuff-3}) 
established {\em local} (in time) regularity for any $T\le T_0$, 
with sufficiently small $T_0$. 
However, a slight variant of the proof yields the estimate \eqref{e:improved}
without any constraint on $T$ (provided it is finite). 
Indeed, it is sufficient to observe that the inequality (58) in the original
proof, that is
\begin{equation}\label{e:will-be-improved}
|U_3|_{H^{\alpha}(\Sigma_{s})} \le K\,
\big(|u|_{H^{\alpha}(\Sigma_{s})}T^{1+2\epsilon} + 
+ |w_0|_{1+\alpha,\Omega_s} + |w_1|_{\alpha,\Omega_s}\big)\,,
\end{equation}
can be actually improved to 
\begin{equation*}
|U_3|_{H^{\alpha}(\Sigma_{s})} \le \epsilon 
|u|_{H^{\alpha}(\Sigma_{s})} + C_{\epsilon} |u|_{L_2(\Sigma_s)} 
+ K \big( |w_0|_{1+\alpha, \Omega_s} + |w_1|_{\alpha,\Omega_s} \big)\,.
\end{equation*}
The above estimate follows as a consequence of the following interpolation 
inequality, 
\begin{equation*} 
\int_0^T |D_t^{\alpha - \epsilon} u |_{L_2(\Gamma_s)} ^2 dt
\leq \epsilon \int_0^T |D_t^{\alpha } u |_{L_2(\Gamma_s)} ^2 dt
+ C_{\epsilon}  |u |_{L_2(\Gamma_s)} ^2 dt\,,
\end{equation*}
which makes it possible to avoid the use of $T^{1/p}$ in order to control the size 
of the principal term in the estimates preceding \eqref{e:will-be-improved}.
}
\end{remark}

\end{document}